\documentclass[10pt]{article}

\usepackage{amsmath,amsfonts,amssymb,amsthm}
\usepackage{a4wide}

\newcommand{\N}{\mathbb{N}}
\newcommand{\R}{\mathbb{R}} 
\newcommand{\rd}{\mathrm{d}}
\newcommand{\e}{\mathrm{e}}
\newcommand{\ini}{\mathrm{I}}

 \usepackage{todonotes}

\usepackage{authblk}
\author{Fabio A. C. C. Chalub}
\affil{Center for Mathematics and Applications (NovaMath), FCT NOVA and Department of Mathematics, FCT NOVA. Faculdade de Ci\^encias e Tecnologia, Universidade Nova de Lisboa, Quinta da Torre, 2829-516 Caparica, Portugal\\ \texttt{facc@fct.unl.pt}}

\title{Solution of the Neutral Kimura equation with two integral constraints}
\date{\today}

\begin{document}

\maketitle

\begin{abstract}
 The Kimura equation is a degenerated partial differential equation of drift-diffusion type used in population genetics. Its solution is required to satisfy not only the equation but a series of conservation laws formulated as integral constraints. In this work, we consider a population of two types evolving without mutation or selection, the so-called neutral evolution. We obtain explicit solutions in terms of Gegenbauer polynomials. To satisfy the integral constraints it is necessary to prove new relations satisfied by the Gegenbauer polynomials. The long-term in time asymptotics is also studied.
 \end{abstract}



\section{Introduction}

The Kimura equation is a partial differential equation of drift-diffusion type introduced in~\cite{Kimura_1954} in the framework of population genetics. The solution of the Kimura Equations, $p(x,t)$, represents the probability density of finding that, at time $t>0$, a fraction $x\in[0,1]$ of individuals carries a given allele, in a population in which two alleles are present. The solution of the neutral case (i.e., when the two alleles do not present any reproductive difference) is known to be written in terms of Gegenbauer polynomials, cf.~\cite{Kimura_1955,Jensen_1974,Aase_1976,Ewens_2004,crow1970introduction}. 

In~\cite{ChalubSouza09b} and~\cite{ChalubSouza14a} the Kimura equation was formally derived as the large population, small time-step limit, limit of the Moran and the Wright-Fisher processes, two well-known finite population Markov processes used in population genetics. See~\cite{ChalubSouza_2017} and references therein for the definition of these two processes.

These two references in the preceding paragraph, by the author and collaborator, generalized the original Kimura equation to allow arbitrary interaction among individuals in the population; see~\eqref{eq:Kimura} below. However, to obtain the correct solution, two additional conservation laws had to be introduced into the model.

The finite population Markov processes referred above are associated with a number of conservation laws to be satisfied by their solutions. Most notably, the solution of the Markov process should satisfy the conservation of probability, but, in fact, the number of conservation laws is equal to the number of different alleles modeled by the associated Markov process. In the present paper, we will consider only the case of two alleles (and therefore, the state of the population will be identified by a single variable, $x\in[0,1]$, that indicates the presence of the focal allele). The correct solution of the continuous models, in the sense of being a continuous approximation of the discrete process, is a measure that not only satisfy the Kimura equation~\eqref{eq:Kimura} but also the two conservation laws given by Eqs.~\eqref{eq:CLs}.

The present work aims to explicitly present the solution of the neutral Kimura equation that satisfies the associated conservation laws. The main difference between this and previous works on the same subject is that the conservation laws are taken into consideration. The construction of the correct solution will require proof of what is, to the best of our knowledge, a new formula involving integrals of Gegenbauer polynomials.  See Eqs.~\eqref{eq:GP_odd}--\eqref{eq:GP_xodd}.

\subsection{Outline}

In Section~\ref{sec:preliminaries} we review some well-known results for the Gegenbauer polynomials and some results from the author and collaborators on the Kimura equation. This section will also fix the notation, as even in the case of Gegenbauer polynomials the notation seems to not be used in consistently\footnote{For example, in~\cite{MorseFeshbach} and in the references cited in the introduction, the $\alpha$-Gegenbauer polynomials are called $T^\alpha_n$. On the other hand, in~\cite{GradshteynRyzhik} they are referred as $C^\alpha_n$, with $T^\alpha_n=C^{\alpha+\frac{1}{2}}_n$. In this work, we use the $C^\alpha_n$ notation, which is also used in large reference websites, such as Wikipedia and Wolfram MathWorld.}. We state the new formulas for $\int_{-1}^1C^\alpha_n(x)\rd x$ and $\int_{-1}^1xC^\alpha_n(x)\rd x$, but the proofs are deferred to the Appendix~\ref{ap:new_formulas}.

In Section~\ref{sec:solution}, we obtain the explicit solution of the Kimura equation that satisfies both conservation laws and derive asymptotic expressions for $t\to\infty$. 

\section{Preliminaries}
\label{sec:preliminaries}

\subsection{Gegenbauer polynomials}

Gegenbauer polynomials $C_n^\alpha:[-1,1]\to\R$, $\Re\alpha>-\frac{1}{2}$, $n\in\N$, are solutions of the Gegenbauer differential equation 
\[
(1-x^2)y''-(2\alpha+1)xy'+n(n+2\alpha)y=0\ .
\]
See~\cite{GradshteynRyzhik}. The set $\{C_n^\alpha\}_{n\in\N}$ is orthogonal, with appropriate weight:
\begin{equation}\label{eq:GP_orthogonality}
\int_{-1}^1C_n^\alpha(x)C_m^\alpha(x)(1-x^2)^{\alpha-\frac{1}{2}}\rd x=\delta_{nm}\frac{\pi 2^{1-2\alpha}\Gamma(n+2\alpha)}{n!(n+\alpha)|\Gamma(\alpha)|^2}\ ,
\end{equation}
where $\delta_{nm}$ is the Kronecer delta. The generating function is given by
\begin{equation}\label{eq:GP_generation}
\frac{1}{\left(1-2xt+t^2\right)^\alpha}=\sum_{n=0}^\infty C_n^\alpha(x)t^n\ .
\end{equation}

We will prove at Appendix~\ref{ap:new_formulas} that for $\alpha\ne1$, 
\begin{align}
\label{eq:GP_odd}
&\int_{-1}^1C_{2n+1}^\alpha(x)\rd x=0\ ,\\
\label{eq:GP_even}
&\int_{-1}^1C_{2n}^\alpha(x)\rd x=\frac{1}{\alpha-1}\binom{2\alpha+2n-2}{2n+1}
\end{align}
and
\begin{align}
\label{eq:GP_xeven}
&\int_{-1}^1 xC_{2n}^\alpha(x)=0\ ,\\
\label{eq:GP_xodd}
&\int_{-1}^1xC_{2n+1}^\alpha(x)\rd x=\frac{2(\alpha+n)}{(\alpha-1)(2n+3)}\binom{2\alpha+2n-2}{2n+1}\ .
\end{align}
The singularity $\alpha=1$ at the right hand side of Eqs.(\ref{eq:GP_odd})--(\ref{eq:GP_xodd}) are removable and the equations can easily be extended to that case.

To the best of our knowledge, Eqs.~(\ref{eq:GP_odd})--(\ref{eq:GP_xodd}) are new. For $\alpha<3/2$, the right hand side of Eqs.~(\ref{eq:GP_odd})--(\ref{eq:GP_xodd}) are identically zero.

\subsection{The Kimura equation}

The Kimura equation is a degenerated partial differential equation of drift-diffusion type (Fokker-Planck), introduced as a model for the evolution of allele frequencies in a given population:
\begin{equation}\label{eq:Kimura}
\partial_tp=\kappa\partial_x^2\left(x(1-x)p\right)-\partial_x\left(x(1-x)\psi(x)p\right)\ .
\end{equation}
In the above equation, $p$ is the probability density that a given allele is present in a fraction $x\in[0,1]$ of individuals in a population, at time $t$, and $\psi:[0,1]\to\R$ is the fitness difference between the focal type and the alternative type, as a function of the presence of the focal type.
 
 Eq.~\eqref{eq:Kimura} has a classical solution $r\in C^\infty(0,1,\R)$, discussed in~\cite{Feller_1951}; see also~\cite{Chalub_Souza:CMS_2009,Danilkina_etal_2018}.
 When $\psi=0$, we call Equation~(\ref{eq:Kimura}) as the \emph{neutral} Kimura equations.

However, the classical solution decay in the limit $t\to\infty$, and therefore cannot be the correct solution from the modeling point of view. Namely, it cannot be a continuous approximation of a given stochastic model, when the number of interacting individuals in a population is very large.
In~\cite{Chalub_Souza:CMS_2009}, solutions of measure type were considered, including Dirac deltas supported on the boundaries of the domain, $x=0$ and $x=1$. More precisely, 
\begin{equation}\label{eq:full_solution}
p(x,t)=a(t)\delta_0(x)+r(x,t)+b(t)\delta_1(x)\ ,
\end{equation}
where $a,b:\R_+\to\R$ are such that two conservation laws derived directly from the discrete process are satisfied: $\partial_t\int_0^1p(x,t)\rd x=\partial_t\int_0^1\varphi(x) p(x,t)\rd x=0$,
where $\varphi$ is the solution of $\varphi''+\psi\varphi'=0$, $\varphi(0)=0$, $\varphi(1)=1$. See~\cite{ChalubSouza14a} for further discussions. In the neutral case $\psi=0$, and therefore, the two conservation laws are
\begin{equation}\label{eq:CLs}
\partial_t\int_0^1p(x,t)\rd x=\partial_t\int_0^1x p(x,t)\rd x=0\ ,
\end{equation}

\section{Solution of the Kimura equation}
\label{sec:solution}


Consider a smooth initial condition $p^\ini$ such that $\mathop{\mathrm{supp}} p^I\subset(0,1)$.
The Gegenbauer polynomials $C_n^{\alpha}(2x-1)$ are solution of $x(1-x)f''+\left(\alpha+\frac{1}{2}\right)(1-2x)f'+n(n+2\alpha)f=0$. It is clear that
\[
r(x,t)=\sum_{n=0}^\infty d_n\e^{-(n+1)(n+2)t}C_n^{3/2}(2x-1)
\]
is a classical solution of the neutral Kimura equation, for any choice of the coefficients $d_n$, $n\in\N$. Let $p$ be given by \eqref{eq:full_solution}, such that conservation laws~\eqref{eq:CLs} are satisfied. Therefore:
\[
a(t)+b(t)+\int_0^1 r(x,t)\rd x=\int_0^1p^\ini(x)\rd x=1\ ,
\]
and
\[
b(t)+\int_0^1 xr(x,t)\rd x=\int_0^1 xp^\ini(x)\rd x\ .
\]
 
 We note that for $\alpha=\frac{3}{2}$
 \[
 \int_0^1 C_n^{3/2}(2x-1)\rd x=\left\{\begin{array}{ll}
 1\ ,&\quad n\in 2\N\ ,\\
 0\ ,&\quad n\in 2\N+1\ ,
 \end{array}\right.
 \]
 and $\int_{0}^1 xC_n^{3/2}(2x-1)\rd x=\frac{1}{2}$ for all $n\in\N$.

The initial condition is given by
$p^\ini(x)=\sum_{n=0}^\infty d_nC_n^{3/2}(2x-1)$. From the orthogonality property of the Gegenbauer polynomials, we find the coefficients
$\left(d_n\right)_{n\in\N}$
\[
d_n=\frac{4(2n+3)}{(n+1)(n+2)}\int_0^1C_n^{3/2}(2x-1)p^\ini(x)x(1-x)\rd x\ .
\]
On the other hand, 
\[
1=\int_0^1p^\ini(x)\rd x=\sum_{n=0}^\infty d_{2n}\ ,
\]
and
\[
\int_0^1xp^\ini(x)\rd x=\frac{1}{2}\sum_{n=0}^\infty d_n\ .
\]

 Finally, the time-dependent fixation probability is given by
 \[
 b(t)=\int xp^\ini(x)\rd x-\int_0^1xr(x,t)\rd x=\frac{1}{2}\sum_{n=0}^\infty d_n\left(1-\e^{-(n+1)(n+2)t}\right)\ ,
 \]
 and time-dependent extinction probability is given by
 \begin{align*}
 a(t)&=1-\int_0^1r(x,t)\rd x-b(t)=\sum_{n\in2\N} d_n-\sum_{n\in 2\N}d_n\e^{-(n+1)(n+2)t}-\frac{1}{2}\sum_{n=0}^\infty d_n\left(1-\e^{-(n+1)(n+2)t}\right)\\
 &=\frac{1}{2}\sum_{n=0}^\infty (-1)^{n}d_n\left(1-\e^{-(n+1)(n+2)t}\right)\ .
 \end{align*}
 
 When $t\to\infty$, 
 \[
  b(t)\approx \int_0^1 xp^\ini(x)\rd x-3\left[\int_0^1p^\ini(x)x(1-x)\rd x\right]\,\e^{-2t}\ ,
 \]
where it was used that $C_0^\alpha(x)=1$. If we assume (in the sense of distributions) that $p^\ini=\delta_{x_0}$, i.e., an initial condition of Dirac-delta type, we conclude that $b(t)\approx x_0-3x_0(1-x_0)\e^{-2t}$, when $t\to\infty$.

 
\appendix

\section{Properties of Gegenbauer polynomials}
\label{ap:new_formulas}

%

Let us define 
\[
f(\alpha,t)=\int_{-1}^1\frac{\rd x}{\left(1-2xt+t^2\right)^\alpha}\ ,\quad \alpha\ne 1\ .
\]
The value of $f$ will be determined using Feynmann's trick for integral calculation\footnote{It is not easy to find classical references to Feynmann's trick. Two online references that discuss this technique in detail are
http://fy.chalmers.se/\~{ }tfkhj/FeynmanIntegration.pdf (Anonymous) and https://courses.grainger.illinois.edu/PHYS487/sp2021/homework/Hwk04Assistance-Feynman\%27sTrick.pdf, by Saavanth Velury. (Both consulted on June 21, 2022.)
}. Differentiating with respect to $\alpha$, and, after one integration by parts, one finds:
\begin{align*}
\partial_\alpha f(\alpha,t)&=-\int_{-1}^1\frac{\log(1-2xt+t^2)}{\left(1-2xt+t^2\right)^\alpha}\rd x\\
&=-\left.\frac{\log(1-2x+t^2)}{(\alpha -1)2t\left(1-2xt+t^2\right)^{\alpha-1}}\right|_{x=-1}^{x=1}-\frac{1}{\alpha-1}\int_{-1}^1\frac{\rd x}{(1-2xt+t^2)^\alpha}\ .
\end{align*}
Therefore, $f$ satisfies 
\[
\partial_\alpha f+\frac{1}{\alpha-1}f=-\frac{1}{(\alpha-1)t}\left[\frac{\log(1-t)}{(1-t)^{2\alpha-2}}-\frac{\log(1+t)}{(1+t)^{2\alpha-2}}\right]\ ,
\]
with initial condition $f(0,t)=2$ for all values of $t$.

The solution of the differential equation is given by
\[
f(\alpha,t)=\frac{1}{\alpha -1}\left[C(t)-\frac{1}{2t}\left(\frac{1}{(1+t)^{2\alpha-2}}-\frac{1}{(1-t)^{2\alpha-2}}\right)\right]\ ,
\]
for a certain $\alpha$-independent function $C$. Using the boundary conditions, we conclude that 
\[
f(\alpha,t)=-\frac{1}{(\alpha-1)2t}\left[\frac{1}{(1+t)^{2\alpha-2}}-\frac{1}{(1-t)^{2\alpha-2}}\right]\ .
\]
From the Taylor expansion
\[
\frac{1}{(1+t)^\beta}=\sum_{n=0}^\infty (-1)^n\binom{\beta+n-1}{n}t^n\ ,
\]
we conclude that
\[
f(\alpha,t)=\frac{1}{\alpha-1}\sum_{n=0}^\infty\binom{2\alpha+2n-2}{2n+1}t^{2n}\ .
\]
and this proves equations~(\ref{eq:GP_odd}) and~(\ref{eq:GP_even}). In particular, if $\alpha<\frac{3}{2}$, $\alpha\ne1$ (including the case of Legendre's polynomials, $\alpha=\frac{1}{2}$), then $\int_{-1}^1C_n^\alpha(x)=0$ for all $n\in\N$. 

Differentiating $f(\alpha,t)$ with respect to $t$:
\[
\partial_t f(\alpha,t)=2\alpha\int_{-1}^1\frac{(x-t)\,\rd x}{(1-2xt+t^2)^{\alpha+1}}=2\alpha\int_{-1}^1\frac{x\,\rd x}{(1-2xt+t^2)^{\alpha+1}}-2\alpha t f(\alpha+1,t)\ ,
\]
we conclude that
\begin{align*}
\int_{-1}^1\frac{x\,\rd x}{(1-2xt+t^2)^\alpha}&=\frac{1}{2(\alpha-1)}\partial_t f(\alpha-1,t)+tf(\alpha,t)\\
&=\frac{1}{2(\alpha-1)(\alpha-2)}\sum_{n=1}^\infty 2n\binom{2\alpha+2n-4}{2n+1}t^{2n-1}
+\frac{1}{\alpha-1}\sum_{n=0}^\infty\binom{2\alpha+2n-2}{2n+1}t^{2n+1}\\
&=\frac{1}{\alpha-1}\left[\sum_{n=0}^\infty\left[ \frac{n+1}{\alpha-2}\binom{2\alpha+2n-2}{2n+3}
+\binom{2\alpha+2n-2}{2n+1}\right]t^{2n+1}\right]\\
&=\frac{1}{\alpha-1}\sum_{n=0}^\infty\frac{2(\alpha+n)}{2n+3}\binom{2\alpha+2n-2}{2n+1}t^{2n+1}\ .
\end{align*}
This proves equations~(\ref{eq:GP_xeven}) and~(\ref{eq:GP_xodd}). All expressions for $f$ can be analytically extended to the case $\alpha=1$.

\section*{Ackowledgements}

This work is funded by national funds through the FCT - Funda\c c\~ao para a Ci\^encia e a Tecnologia, I.P., under the scope of the projects UIDB/00297/2020 and UIDP/00297/2020 (Center for Mathematics and Applications).


\begin{thebibliography}{10}

\bibitem{Kimura_1954}
Motoo Kimura.
\newblock {Process Leading To Quasi-Fixation Of Genes In Natural Populations
  Due To Random Fluctuation Of Selection Intensities}.
\newblock {\em Genetics}, 39(3):280--295, 05 1954.

\bibitem{Kimura_1955}
Motoo Kimura.
\newblock Solution of a process of random genetic drift with a continuous
  model.
\newblock {\em Proceedings of the National Academy of Sciences of the United
  States of America}, 41(3):144--150, 1955.

\bibitem{Jensen_1974}
Louis Jensen.
\newblock Solving a singular diffusion equation occurring in population
  genetics.
\newblock {\em J. Appl. Probability}, 11:1--15, 1974.

\bibitem{Aase_1976}
Knut Aase.
\newblock A note on a singular diffusion equation in population genetics.
\newblock {\em J. Appl. Probability}, 13(1):1--8, 1976.

\bibitem{Ewens_2004}
Warren~J. {Ewens}.
\newblock {\em {Mathematical population genetics. I: Theoretical introduction.
  2nd ed.}}
\newblock New York, NY: Springer, 2nd ed. edition, 2004.

\bibitem{crow1970introduction}
James~Franklin Crow, Motoo Kimura, et~al.
\newblock {\em An introduction to population genetics theory.}
\newblock New York, Evanston and London: Harper \& Row, Publishers, 1970.

\bibitem{ChalubSouza09b}
Fabio A. C.~C. Chalub and Max~O. Souza.
\newblock From discrete to continuous evolution models: a unifying approach to
  drift-diffusion and replicator dynamics.
\newblock {\em Theor. Pop. Biol.}, 76(4):268--277, 2009.

\bibitem{ChalubSouza14a}
Fabio A. C.~C. Chalub and Max~O. Souza.
\newblock {The frequency-dependent Wright-Fisher model: diffusive and
  non-diffusive approximations}.
\newblock {\em J. Math. Biol.}, 68(5):1089--1133, 2014.

\bibitem{ChalubSouza_2017}
Fabio A. C.~C. {Chalub} and Max~O. {Souza}.
\newblock {On the stochastic evolution of finite populations.}
\newblock {\em {J. Math. Biol.}}, 75(6-7):1735--1774, 2017.

\bibitem{MorseFeshbach}
Philip~M. Morse and Herman Feshbach.
\newblock {\em Methods of theoretical physics. 2 volumes}.
\newblock McGraw-Hill Book Co., Inc., New York-Toronto-London, 1953.

\bibitem{GradshteynRyzhik}
I.~S. Gradshteyn and I.~M. Ryzhik.
\newblock {\em Table of integrals, series, and products}.
\newblock Academic Press [Harcourt Brace Jovanovich, Publishers], New
  York-London-Toronto, Ont., 1980.
\newblock Corrected and enlarged edition edited by Alan Jeffrey, Incorporating
  the fourth edition edited by Yu. V. Geronimus [Yu. V. Geronimus] and M. Yu.
  Tseytlin [M. Yu. Tse\u{\i}tlin], Translated from the Russian.

\bibitem{Feller_1951}
William Feller.
\newblock Diffusion processes in genetics.
\newblock In {\em Proceedings of the {S}econd {B}erkeley {S}ymposium on
  {M}athematical {S}tatistics and {P}robability, 1950}, pages 227--246.
  University of California Press, Berkeley-Los Angeles, Calif., 1951.

\bibitem{Chalub_Souza:CMS_2009}
Fabio A. C.~C. Chalub and Max~O. Souza.
\newblock A non-standard evolution problem arising in population genetics.
\newblock {\em Commun. Math. Sci.}, 7(2):489--502, 2009.

\bibitem{Danilkina_etal_2018}
Olga Danilkina, Max~O. Souza, and Fabio A. C.~C. Chalub.
\newblock Conservative parabolic problems: nondegenerated theory and
  degenerated examples from population dynamics.
\newblock {\em Math. Methods Appl. Sci.}, 41(12):4391--4406, 2018.

\end{thebibliography}

\end{document}